# 2-Rainbow domination number of circulant graphs $C(n; \{1,4\})$

[a] **Ramy Shaheen,** [b] **Suhail Mahfud** and [c] **Mohammed Fahed Adrah**

[a,b,c] Department of Mathematics, Faculty of science Tishreen University, Lattakia, Syria.

[a] Aljazeera Private University, Damascus, Syria.

[a] E-mail: shaheenramy2010@hotmail.com, [b] E-mail: mahfudsuhail@gmail.com, [c] E-mail: mohadrah7@gmail.com.

## Abstract

Let k be a positive integer, a k-rainbow domination function ($kRDF$) of a graph G is a function $f$ from $V(G)$ to the set of all subsets of $\{1,2,\ldots,k\}$ such that every vertex $v \in V(G)$ with $f(v) = \emptyset$ satisfies $\bigcup_{u \in N(v)} f(u) = \{1, 2, \ldots, k\}$. The weight of a $kRDF$ is defined as the value $w(f) = \sum_{v \in V(G)} |f(v)|$. The k-rainbow domination number of $G$, which is denoted by $\gamma_{rk}(G)$, is the minimum weight of all $kRDFs$ of $G$. In this paper we were able to give the exact value of the 2-rainbow domination number of circulant graphs $C(n; \{1, 4\})$, which is $\gamma_{r2}(C(n; \{1, 4\})) = \left\lceil \frac{n}{3} \right\rceil + \alpha$, where $\alpha = 0$ for $n \equiv 0 \pmod{6}$, $\alpha = 1$ for $n \equiv 1, 2, 3, 5 \pmod{6}$ and $\alpha = 2$ for $n \equiv 4 \pmod{6}$.

## 1. Introduction

Let $G(V, E)$ be a simple, finite and undirected graph where $V(G)$ is the vertex set, and $E(G)$ is the edge set. The open neighborhood of a vertex $v$ is $N(v) = \{u \in V(G): uv \in E(G)\}$, and the closed neighborhood is $N[v] = N(v) \cup \{v\}$.

Domination is one of the most important concepts in graph theory, because it has varied applications in many fields such as optimization, operation research, biological networks, communication networks, coding theory, etc. A set $D$ of vertices of a graph $G$ is said to be a dominating set if every vertex in $V - D$ is adjacent to a vertex in $D$. A dominating set $D$ is said to be a minimal dominating set if no proper subset of $D$ is a dominating set. The minimum cardinality of a dominating set of a graph $G$ is called the domination number of $G$ and is denoted by $\gamma(G)$.

The problem of domination in real life is affected by the factors of the issue, which led to the existence of more than 75 different types of domination.

Therefore Brešar et al. in 2005 [1] introduced the k-rainbow domination concept, and when $k = 2$ a 2-rainbow domination function $f$ (2RDF) of a graph G is a function $f$ from $V(G)$ to the set of all subsets of $\{1, 2\}$, such that for any vertex $v$ with $f(v) = \emptyset$ we have $\bigcup_{u \in N(v)} f(u) = \{1, 2\}$. The weight of the function $f$ denoted by $w(f)$ and defined as $w(f) = \sum_{v \in V(G)} |f(v)|$. The 2-rainbow domination number of $G$ (denoted by $\gamma_{r2}(G)$) is the minimum weight of all $2RDFs$ of $G$.

The importance of rainbow domination comes from the fact that most real-life problems require many conditions to say that the problem has been dominated and the issue is solved.

And we will present the following application as an example to illustrate the use of the 2-rainbow domination function in real life.

**Application 1.1.**

Consider the map in the figure 1 is a map of a new suburban. We want to distribute a number of traffic police and traffic lights on the crossroads of the main streets at the lowest possible cost with the possibility of having a policeman and a traffic light at the same intersection, Where the vacant places must have at least one policeman and one traffic light at the crossroads adjacent to them, so if there is any problem the policeman can move between the empty areas adjacent to its own to solve it, While the sign maintains the traffic order in that area. So, first we transform the map into a graph by replacing every intersection with a vertex and every street with an edge, and second, we apply a 2-rainbow domination function on that graph (the numbers 0, 1, 2, and 3 represent the crossroads that is empty, has a policeman, has a traffic sign and has both a traffic sign and a policeman respectively (see Fig 1)).

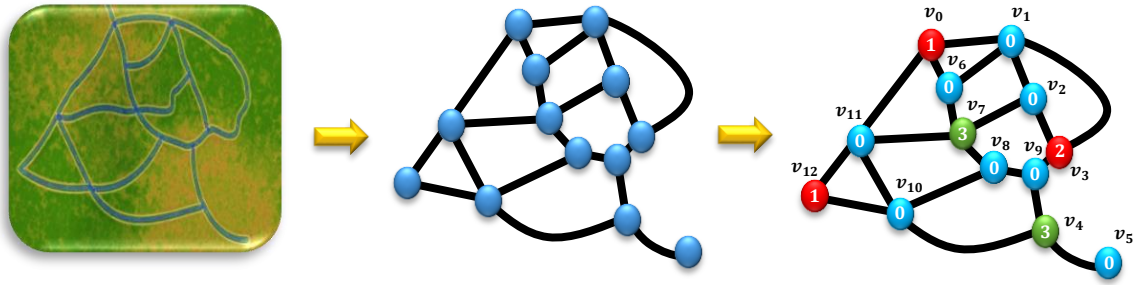

Fig 1

In [2], Brešar and Sumenjak showed that the problem of deciding if a graph has a 2-rainbow dominating function of a given weight is NP-complete. Chunling et al. in [3], were the first to present the proof technique used in this paper, they found the 2-rainbow domination number of generalized Peterson $P(n, 2)$ as the following:

**Theorem 1.**

$$\gamma_{r2}(P(n, 2)) = \begin{cases} \left\lceil \dfrac{4n}{5} \right\rceil: & n \equiv 0, 3, 4, 9 \ (mod \ 10), \\ \left\lceil \dfrac{4n}{5} \right\rceil + 1: & n \equiv 1, 2, 5, 6, 7, 8 \ (mod \ 10). \end{cases}$$

In 2021, Ervis et al. [4], obtained new results on 2-rainbow domination number of generalized Petersen graphs $P(5k, k)$ as the following:

**Theorem 2.** Let k > 3, Then

$$\gamma_{r2}\left(P(5k,k)\right) = \begin{cases} 4k: & k \equiv 2, 8 \ (mod\ 10), \\ 4k+1: & k \equiv 5, 9 \ (mod\ 10). \end{cases}$$

$$4k+1 \leq \gamma_{r2}\left(P(5k,k)\right) \leq \begin{cases} 4k+2: & k \equiv 1, 6, 7 \ (mod\ 10), \\ 4k+3: & k \equiv 0, 3, 4 \ (mod\ 10). \end{cases}$$

The cases when k ≤ 3 are summarized below:

$$\gamma_{r2}\left(P(5,k)\right) = 5, \ \gamma_{r2}\left(P(10,2)\right) = 10, 13 \leq \gamma_{r2}\left(P(15,3)\right) \leq 14.$$

This concept has attracted attention of several authors. See, for example [2-13]. In 2011 Ali et al. [12], were able to prove the following result:

**Theorem 3.** Let G be a K-regular graph, then $\gamma_{r2}(G) \geq \left\lceil \frac{2n}{K+2} \right\rceil$.

For $1 \leq j < k \leq \lfloor \frac{n}{2} \rfloor$, the Circulant graph is a graph with $n$ vertices defined as:

$$V\big(C(n;\{j,k\})\big) = V(G) = \{v_i: 0 \leq i \leq n-1\},$$

$$E\big(C(n;\{j,k\})\big) = E(G) = \{v_i v_{i\pm j}, v_i v_{i\pm k}: 0 \leq i \leq n-1\}.$$

where all the indices $i \pm j$, $i \pm k$ are reduced modulo $n$.

In 2015, Fu et al. [14], found the 2-rainbow domination number of circulant graphs $C(n;\{1,3\})$ as the following:

**Theorem 4.** For $n \geq 7$, Let $m = \lfloor \frac{n}{5} \rfloor$ and $n = 5m + \alpha$ then

$$\gamma_{r2}\left(C(n;\{1,3\})\right) = \begin{cases} 2m: & \alpha = 0, \\ 2m+1: & \alpha = 1, 2, \\ 2m+2: & \alpha = 3, 4. \end{cases}$$

Domination and its variation of the circulant graphs have been studied a lot in the past years see [14-16], and from the definition of 2-rainbow domination, various sorts of dominating parameters were provided, such as (Maximal 2-rainbow domination [17], Restrained 2-rainbow domination [18, 19], Total 2-rainbow domination [20, 21], Independent 2-rainbow domination [22, 23]).

## 2. 2-rainbow domination number of $C(n;\{1,4\})$

**Lemma 2.1.** For $n \geq 9$ then

$$\gamma_{r2}\left(C(n;\{1,4\})\right) \leq \begin{cases} \left\lceil \frac{n}{3} \right\rceil: & n \equiv 0 \ (mod\ 6), \\ \left\lceil \frac{n}{3} \right\rceil + 1: & n \equiv 1, 2, 3, 5 \ (mod\ 6), \\ \left\lceil \frac{n}{3} \right\rceil + 2: & n \equiv 4 \ (mod\ 6). \end{cases}$$

**Proof.** To show the upper bound we define the function $f$ as follows:
$$f(v_0, v_1, \ldots, v_{n-1}) = (f(v_0)f(v_1)\ldots f(v_{n-1})).$$

where $0, 1, 2, 3$ stand for $\emptyset, \{1\}, \{2\}, \{1,2\}$ respectively, for example we will apply $f$ on $C(19; \{1,4\})$ as follows:
$$f(v_0, v_1, \ldots, v_{18}) = (100200 \ldots 1002201) \text{ (See Fig 2.1)}.$$

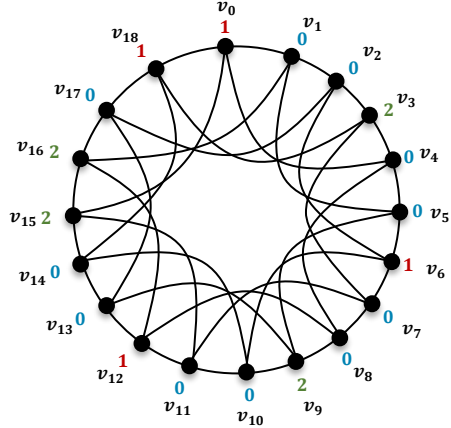

Fig 2.1 $C(19; \{1,4\})$

Now we have:
$$f(v_0, v_1, \ldots, v_{n-1}) = \begin{cases} (100200 \ldots 100200): & n \equiv 0 \pmod 6, \\ (100200 \ldots 1002001002201): & n \equiv 1 \pmod 6, \\ (100200 \ldots 10020010020210): & n \equiv 2 \pmod 6, \\ (300200100200 \ldots 100200100): & n \equiv 3 \pmod 6, \\ (100200 \ldots 1002001212): & n \equiv 4 \pmod 6, \\ (100200 \ldots 10020010220): & n \equiv 5 \pmod 6. \end{cases}$$

And since $f$ is 2RDF of $C(n; \{1,4\})$ with

$$w(f) = \begin{cases} \left\lceil \frac{n}{3} \right\rceil: & n \equiv 0 \pmod 6, \\ \left\lceil \frac{n}{3} \right\rceil + 1: & n \equiv 1, 2, 3, 5 \pmod 6, \\ \left\lceil \frac{n}{3} \right\rceil + 2: & n \equiv 4 \pmod 6. \end{cases}$$

Then we have

$$\gamma_{r2}(C(n; \{1,4\})) \leq \begin{cases} \left\lceil \frac{n}{3} \right\rceil: & n \equiv 0 \pmod 6, \\ \left\lceil \frac{n}{3} \right\rceil + 1: & n \equiv 1, 2, 3, 5 \pmod 6, \\ \left\lceil \frac{n}{3} \right\rceil + 2: & n \equiv 4 \pmod 6. \end{cases} \qquad \square$$

**Proposition 2.2.** We will define the following sets to use in the proof:

$V_0 = \{v \in V(C(n; \{1,4\})): f(v) = \emptyset\}$,

$V_1 = \{v \in V(C(n; \{1,4\})): f(v) \in \{\{1\},\{2\}\}\}$,

$V_2 = \{v \in V(C(n; \{1,4\})): f(v) = \{\{1,2\}\}\}$,

$V_{ij} = \{v \in V_0: i$ denotes to the number of vertices adjacent to $v$ from $V_1$ and $j$ denotes to the number of vertices adjacent to $v$ from $V_2\}$,

$E_1 = \{uv \in E(C(n; \{1,4\})): u, v \in V_1\}$,

$E_2 = \{uv \in E(C(n; \{1,4\})): u, v \in V_2\}$,

$E_{12} = \{uv \in E(C(n; \{1,4\})): u \in V_1, v \in V_2\}$.

Then $V(C(n; \{1,4\})) = V_0 \cup V_1 \cup V_2$, where $V_i \cap V_j = \emptyset$, $i, j = 0, 1, 2$ and $i \neq j$.

Let $T = \{V_{01}, V_{02}, V_{03}, V_{04}, V_{11}, V_{12}, V_{13}, V_{20}, V_{21}, V_{22}, V_{30}, V_{31}, V_{40}\}$.

Where the collection $T$ is pairwise disjoint, and $V_0 = \bigcup_{S \in T} S$.

The number of edges that connects between $V_0$ and $V_1$ is $4|V_1| - 2|E_1| - |E_{12}|$ and this number also can be presented by $|V_{11}| + |V_{12}| + |V_{13}| + 2|V_{20}| + 2|V_{21}| + 2|V_{22}| + 3|V_{30}| + 3|V_{31}| + 4|V_{40}|$.

Therefore, we have

$$4|V_1| - 2|E_1| - |E_{12}| = |V_{11}| + |V_{12}| + |V_{13}| + 2|V_{20}| + 2|V_{21}| + 2|V_{22}| + 3|V_{30}| + 3|V_{31}| + 4|V_{40}|. \quad (1)$$

And in the same way we have

$$4|V_2| - 2|E_2| - |E_{12}| = |V_{01}| + 2|V_{02}| + 3|V_{03}| + 4|V_{04}| + |V_{11}| + 2|V_{12}| + 3|V_{13}| + |V_{21}| + 2|V_{22}| + |V_{31}|. \quad (2)$$

Furthermore by (1) +2× (2) we obtain

$$4|V_1| + 8|V_2| = 3|V_{11}| + 5|V_{12}| + 7|V_{13}| + 2|V_{20}| + 4|V_{21}| + 6|V_{22}| + 3|V_{30}| + 5|V_{31}| + 4|V_{40}| + 2|V_{01}| + 4|V_{02}| + 6|V_{03}| + 8|V_{04}| + 3|E_{12}| + 2|E_1| + 4|E_2|.$$

And since $n = |V_0| + |V_1| + |V_2|$ then $2|V_0| = 2n - 2|V_1| - 2|V_2|$, which leads to

$$4|V_1| + 8|V_2| = 2n - 2|V_1| - 2|V_2| + |V_{11}| + 3|V_{12}| + 5|V_{13}| + 2|V_{21}| + 4|V_{22}| +$$
$$|V_{30}| + 3|V_{31}| + 2|V_{40}| + 2|V_{02}| + 4|V_{03}| + 6|V_{04}| + 3|E_{12}| + 2|E_1| +$$
$$4|E_2| + 2|V_2| - 2|V_2|.$$
$$6|V_1| + 12|V_2| = 2n + |V_{11}| + 3|V_{12}| + 5|V_{13}| + 2|V_{21}| + 4|V_{22}| + |V_{30}| + 3|V_{31}| +$$
$$2|V_{40}| + 2|V_{02}| + 4|V_{03}| + 6|V_{04}| + 3|E_{12}| + 2|E_1| + 4|E_2| + 2|V_2|.$$
$$6w(f) = 2n + |V_{11}| + 3|V_{12}| + 5|V_{13}| + 2|V_{21}| + 4|V_{22}| + |V_{30}| + 3|V_{31}| + 2|V_{40}| +$$
$$2|V_{02}| + 4|V_{03}| + 6|V_{04}| + 3|E_{12}| + 2|E_1| + 4|E_2| + 2|V_2|.$$

Thus, we have Lemma 2.3.

**Lemma 2.3.** $6w(f) = 2n + \beta$ where $\beta = |V_{11}| + 3|V_{12}| + 5|V_{13}| + 2|V_{21}| + 4|V_{22}| + |V_{30}| + 3|V_{31}| + 2|V_{40}| + 2|V_{02}| + 4|V_{03}| + 6|V_{04}| + 3|E_{12}| + 2|E_1| + 4|E_2| + 2|V_2|.$ □

**Lemma 2.4.** $\beta \geq 6$ when $n \equiv 1, 2, 3, 5 \pmod{6}$ and $\beta \geq 12$ when $n \equiv 4 \pmod{6}$

**Proof.** Suppose the contrary that $\beta \leq 5$ when $n \equiv 1, 2, 3, 5 \pmod{6}$ and $\beta \leq 11$ when $n \equiv 4 \pmod{6}$. (We will define the symbol (~) and use it instead of the phrase (and all other sets in $\beta$ are empty) to avoid repetition).

**Case 1:** $v_0 \in V_2$ and $|N(v_0) \cap V_1| = 3$.

This case only appears when $n \equiv 4 \pmod{6}$, it makes $|V_2| = 1$, $|E_{12}| = 3$, and $|V_{11}| = |V_{12}| = |V_{13}| = |V_{21}| = |V_{22}| = |V_{30}| = |V_{31}| = |V_{40}| = |V_{02}| = |V_{03}| = |V_{04}| = |E_1| = |E_2| = 0$.

Let $v_1, v_{n-1}, v_4 \in V_1$, It makes $v_{n-2}, v_{n-3}, v_{n-4}, v_{n-5}, v_2, v_3, v_5, v_8 \in V_0$ (otherwise either $|E_1| \neq 0$ or $|E_2| \neq 0$), $v_2$ makes $v_6 \in V_1$ (otherwise a contradiction with 2RDF definition). Thus, we have $v_5 \in V_{30} \cup V_{40}$ (see Fig 2.2 where Gray dot, Black dot and Xdot (a white dot with an X inside) stand for the vertex of $V_0, V_1, V_2$ respectively), a contradiction.

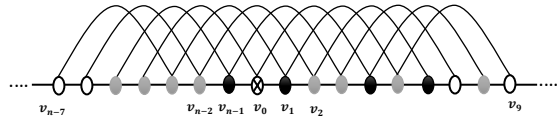

Fig 2.2: $|N(v_0) \cap V_1| = 3$.

**Case 2:** $v_0 \in V_2$ and $|N(v_0) \cap V_1| = 2$.

Also this case only appears when $n \equiv 4 \pmod{6}$, it makes $|V_{13}| = |V_{22}| = |V_{03}| = |V_{04}| = |E_2| = 0$, $|V_{11}| + 3|V_{12}| + 2|V_{21}| + |V_{30}| + 2|V_{02}| + 3|E_{12}| + 2|E_1| + 2|V_2| + 3|V_{31}| + 2|V_{40}| \leq 11$ and $|E_2| < 2$. Now, let $v_1, v_4 \in V_1$ then $v_{n-1} \in V_{01} \cup V_{21} \cup V_{31} \cup V_{11}$ (Otherwise $\beta > 11$), and we will consider all these cases as the following:

- $v_{n-1} \in V_{21}$, it makes $v_2, v_3 \in V_0$ and $v_{n-4} \in V_{11}$ (Otherwise $\beta > 11$). Now, $|V_{11}| = |V_2| = |V_{21}| = 1$, $|E_{12}| = 2\sim$, $v_{n-6}, v_{n-7} \in V_0$ (otherwise either $|E_1| \neq 0$ or $|V_{30}| \neq 0$), but $v_{n-7}$ forces $v_{n-11} \in V_2$ (see Fig 2.3 (1)), a contradiction.

- $v_{n-1} \in V_{11}$ we face three cases:

  ○ $v_3 \in V_1$, it makes $v_3 v_4 \in E_1$, $|V_{11}| = |V_2| = |E_1| = 1$, $|E_{12}| = 2\sim$, moreover $v_{n-2}, v_{n-3}, v_2 \in V_0$ (otherwise $E_1 > 1$ or $|V_{30}|, |V_{40}| \neq 0$), $v_{n-2}$ forces $v_{n-6} \in V_2$ (see Fig 2.3 (2)), a contradiction.

  ○ $v_{n-2} \in V_1$, it makes $v_{n-3}, v_2 \in V_0$ (otherwise $|E_1| > 1$), $v_{n-5}$ forces $v_{n-9}, v_{n-6} \in V_1$ (otherwise a contradiction with 2RDF definition). Thus, we have $v_{n-6} v_{n-2} \in E_1$. Now, $|V_{11}| = |V_2| = |E_1| = 1$, $|E_{12}| = 2\sim$, so $v_{n-7}, v_{n-8}, v_{n-10} \in V_0$ (otherwise $|E_1| > 1$), $v_{n-7}$ forces $v_{n-11} \in V_1$, which implies that $v_{n-10} \in V_{30} \cup V_{40}$ (see Fig 2.3 (3)), a contradiction.

  ○ $v_{n-5} \in V_1$, $v_{n-3} \in V_0$ (otherwise $v_1 v_{n-3} \in E_1$ and $v_{n-4} \in V_{21} \cup V_{31}$ which makes $\beta > 11$). $v_{n-2}$ forces $v_{n-6}, v_2 \in V_1$ (otherwise a contradiction with 2RDF definition), that leads to $v_1 v_2, v_{n-6} v_{n-5} \in E_1$ (see Fig 2.3 (4)), a contradiction.

- $v_{n-1} \in V_{01}$, if $v_{n-3} \in V_1$ then $v_1 v_{n-3} \in E_1$, $v_{n-4} \in V_{11}$, $|V_{11}| = |V_2| = |E_1| = 1$, $|E_{12}| = 2\sim$. Now, $v_2$ forces $v_6 \in V_1$ and that makes $v_5 \in V_{30} \cup V_{40}$ (see Fig 2.3 (5)), a contradiction. Therefore $v_{n-3} \in V_0$ and it forces either $v_{n-7} \in V_2$ then $|V_2| = |E_{12}| = 2$, $|V_{11}| = 1\sim$, and that makes $v_{n-6}, v_{n-8} \in V_0$, so $v_{n-2}$ forces $v_2 \in V_2$ (see Fig 2.3 (6)), a contradiction, or $v_{n-7} \in V_1$ then $v_{n-2}$ forces $v_2, v_{n-6} \in V_1$, $v_1 v_2, v_{n-6} v_{n-7} \in E_1$ and $2|V_2| + 3|E_{12}| + 2|E_1| > 11$ (see Fig 2.3 (7)), a contradiction.

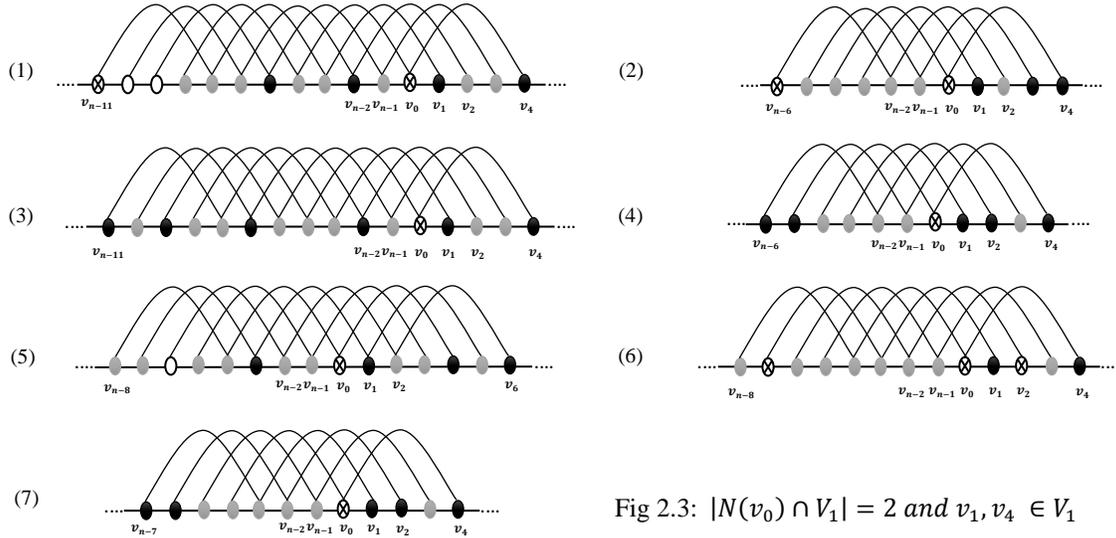

Fig 2.3: $|N(v_0) \cap V_1| = 2$ and $v_1, v_4 \in V_1$

**Case 3:** $v_0 \in V_2$ and $|N(v_0) \cap V_1| = 1$.

**Case 3.1:** $n \equiv 1, 2, 3, 5 \pmod 6$, It follows that $|V_2| = |E_{12}| = 1\sim$. Now, if $v_1 \in V_1$ that makes $v_{n-1}, v_{n-4}, v_4 \in V_{01}$ then $v_3$ forces $v_7 \in V_2$ (see Fig 2.4), a contradiction.

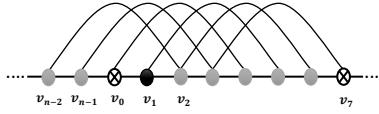

Fig 2.4: $|N(v_0) \cap V_1| = 1, n \equiv 1, 2, 3, 5 \pmod 6$.

**Case 3.2:** $n \equiv 4 \pmod 6$, In this case $|N(v_0) \cap V_2| \leq 1$ otherwise $\beta > 11$.

**Case 3.2.1:** $|N(v_0) \cap V_2| = 1$, then we have $|V_2| = 2, |E_{12}| = |E_2| = 1\sim$. Now, if $v_1 \in V_2$ then $v_{n-1} \in V_1$ (otherwise $2|V_2| + 3|E_{12}| + 4|E_2| + |V_{11}| > 11$) and $v_{n-8}$ forces $v_{n-12} \in V_2$ (see Fig 2.5 (1)), a contradiction, and if $v_4 \in V_2$ then $v_{n-4} \in V_1$ (otherwise $2|V_2| + 3|E_{12}| + 4|E_2| + |V_{11}| > 11$) and $v_2$ forces $v_6 \in V_2$ (see Fig 2.5 (2)), a contradiction.

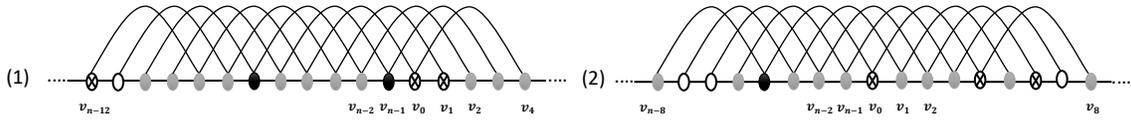

Fig 2.5: $|N(v_0) \cap V_1| = 1, n \equiv 4 \pmod 6, v_1 \in V_2$.

**Case 3.2.2:** $|N(v_0) \cap V_2| = 0$, let $v_1 \in V_1$ then we have $v_{n-4} \in V_{01} \cup V_{11} \cup V_{21} \cup V_{31} \cup V_{02} \cup V_{12}$ and will discuss these cases as follows:

• $v_{n-4} \in V_{12}$ then $v_{n-3} \in V_0$, $v_{n-8} \in V_2$ and $v_{n-5} \in V_1$, and that leads to $v_{n-1} \in V_{11}$, $|V_2| = 2, |E_{12}| = |V_{12}| = |V_{11}| = 1\sim$. $v_{n-2}$ forces $v_{n-6}, v_2 \in V_1$ which makes $|E_1| \neq 0$ (see Fig 2.6 (1)), a contradiction.

• $v_{n-4} \in V_{02}$ let $v_{n-5} \in V_2$, it forces $v_{n-1} \in V_{02}$, that leads to $|V_2| = |V_{02}| = 2, |E_{12}| = 1\sim$. Now, $v_{n-2}$ makes $v_{n-6}, v_2 \in V_1$, thus $|E_{12}| > 1$ (see Fig 2.6 (2)), a contradiction.

• $v_{n-4} \in V_{31}$ then $v_{n-3} v_1 \in E_1$, $v_{n-1} \in V_{11}$, $|V_2| = |E_{12}| = |V_{31}| = |V_{11}| = |E_1| = 1\sim$, so $v_{n-6}, v_2 \in V_0$, $v_{n-2}$ leads to a contradiction with the 2RDF definition (see Fig 2.6 (3)).

• $v_{n-4} \in V_{21}$, let $v_{n-8} \in V_0$ then $v_{n-3} v_1 \in E_1$ and $v_{n-1} \in V_{21} \cup V_{11}$, so if $v_{n-1} \in V_{21}$ then $v_3 \in V_1$, $|V_2| = |E_{12}| = |E_1| = 1$, $|V_{21}| = 2\sim$, but $v_4 \in V_{11} \cup V_{21} \cup V_{31}$ a contradiction, and if $v_{n-1} \in V_{11}$, then $|V_2| = |E_{12}| = |V_{21}| = |E_1| = 1$, $|V_{11}| + |V_{30}| \leq 2\sim$, that leads to $v_2 \in V_0$ and $v_3$ forces $v_7 \in V_2$ (see Fig 2.6 (4)), a contradiction.

• $v_{n-4} \in V_{11}$, let $v_{n-3} \in V_1$ then $v_1 v_{n-3} \in E_1$. Now we will discuss each of the following vertices in all of the possible cases:

  ○ If $v_{n-2} \in V_2$ then $2|V_2| + 3|E_{12}| + 2|E_1| + |V_{11}| > 11$, a contradiction. And if $v_{n-2} \in V_1$ then $v_{n-2} v_{n-3} \in E_1$, $v_{n-1} \in V_{11}$, $|V_2| = |E_{12}| = 1, |E_1| = |V_{11}| = 2\sim$, thus $v_2 \in V_0$ and $v_3$ forces $v_7 \in V_2$ (see Fig 2.6 (5)), a contradiction. Therefore $v_{n-2} \in V_0$.

- If $v_{n-6} \in V_2$ then $v_{n-2} \in V_{11}$, $v_{n-7}, v_{n-10} \in V_0$, $|V_2| = |V_{11}| = 2$, $|E_{12}| = |E_1| = 1\sim$, but $v_{n-7} \in V_{11} \cup V_{21}$ (see Fig 2.6 (6)), a contradiction. If $v_{n-6} \in V_1$ then $v_{n-7} \in V_0$, and in a same way we see that $v_2, v_3 \in V_0$, and $v_3$ forces $v_7 \in V_2$, moreover $v_2$ makes $v_6 \in V_1$, so $2|V_2| + 3|E_{12}| + 2|E_1| + |V_{11}| > 11$ (see Fig 2.6 (7)), a contradiction, Therefore $v_{n-6} \in V_0$ moreover $v_{n-5}$ forces $v_{n-9} \in V_2$ and $v_{n-2}$ makes $v_6 \in V_1$, which implies that $2|V_2| + 3|E_{12}| + 2|E_1| + |V_{11}| > 11$ (see Fig 2.6 (8)), a contradiction.

- $v_{n-4} \in V_{01}$ then $v_{n-1} \in V_{01} \cup V_{11} \cup V_{21} \cup V_{02} \cup V_{12}$ and we will discuss these cases as follows:

  - $v_{n-1} \in V_{12}$ then $|V_2| = 2$, $|E_{12}| = |V_{12}| = 1$, $|V_{11}| + |V_{30}| \leq 1\sim$. Now, if $v_3 \in V_2$ then $v_2 \in V_{21} \cup V_{31}$ (see Fig 2.6 (9)), a contradiction.

  - $v_{n-1} \in V_{02}$, let $v_3 \in V_2$ then $v_2 \in V_0$ and $v_{n-2}$ forces $v_{n-6} \in V_2$, which makes $|V_2| = 3$, $|E_{12}| = |V_{02}| = 1\sim$, but $v_2 \in V_{11} \cup V_{21}$ (see Fig 2.6 (10)), a contradiction.

  - $v_{n-1} \in V_{21}$, if $v_{n-6} \in V_1$ then $v_{n-2} v_{n-6} \in E_1$ and that makes $v_2 \in V_{30}$, $v_4 \in V_{11}$, $|V_2| = |E_{12}| = |E_1| = |V_{30}| = |V_{21}| = |V_{11}| = 1\sim$, therefore $v_6 \in V_0$ and it forces $v_{10} \in V_2$ (see Fig 2.6 (11)), a contradiction. Therefore $v_{n-6} \in V_0$ and $v_{n-5}$ forces $v_{n-9} \in V_2$ which makes $v_2 \in V_{30}$, $v_4 \in V_{11}$, $|V_2| = 2$, $|E_{12}| = |V_{30}| = |V_{21}| = |V_{11}| = 1\sim$, $v_7 \in V_0$ and $v_6$ forces $v_{10} \in V_2$ (see Fig 2.6 (12)), a contradiction.

  - $v_{n-1} \in V_{11}$, let $v_3 \in V_1$, if $v_2 \in V_1$ then $v_1 v_2, v_2 v_3 \in E_1$ and $|V_2| = |E_{12}| = 1$, $|E_1| = 2$, $|V_{11}| + |V_{30}| \leq 2\sim$, but $v_{n-2}, v_{n-3}$ force $v_{n-6}, v_{n-7} \in V_1$ and that makes $v_{n-6} v_{n-7} \in E_1$ (see Fig 2.6 (13)), a contradiction. Therefore $v_2 \in V_0$, $v_{n-2}$ forces $v_{n-6} \in V_2$ and $v_{n-3}$ makes $v_{n-7} \in V_1$ and that leads to $|V_2| = |E_{12}| = 2$, $|V_{11}| = 1\sim$, moreover $v_{n-8}, v_{n-9}, v_{n-10}, v_{n-11} \in V_0$, but $v_{n-9}$ force $v_{n-13} \in V_2$ (see Fig 2.6 (14)), a contradiction.

  - $v_{n-1} \in V_{01}$, let $v_2 \in V_1$ then $v_1 v_2 \in E_1$ and $v_{n-2}$ forces $v_{n-6} \in V_1$ then $v_{n-3}, v_{n-5}$ force $v_{n-7}, v_{n-9} \in V_1$, $v_{n-6} v_{n-7} \in E_1$ and that leads to $v_{n-10}, v_{n-11} \in V_0$, then we face the following cases:

    - $v_{n-12} \in V_1$ then $v_{n-8} \in V_{30}$, $|V_2| = |E_{12}| = 1$, $|E_1| = 2$, $|V_{11}| + |V_{30}| \leq 2\sim$ so $v_{n-13}, v_{n-14} \in V_0$ and $v_{n-14}$ forces $v_{n-15}, v_{n-18} \in V_1$ and that makes $|V_{30}| = 2$, $|V_{11}| = 0$, and by continuing in this way we have for $4 < j < n - 7$

    $$\begin{cases} v_j \in V_0: j = 6i, 6i + 2, 6i + 3, 6i + 5, \\ v_j \in V_1: j = 6i + 1, 6i + 4. \end{cases}$$

    But $v_5$ leads to a contradiction with the 2RDF definition (see Fig 2.6 (15)).

- $v_{n-12} \in V_0$, $v_{n-11}$ forces $v_{n-15} \in V_1$, $v_{n-12}$ makes $v_{n-13}, v_{n-16} \in V_1$ then $v_{n-9}v_{n-13}, v_{n-15}v_{n-16} \in E_1$ and $2|V_2| + 3|E_{12}| + 2|E_1| > 11$ (see Fig 2.6 (16)), a contradiction.

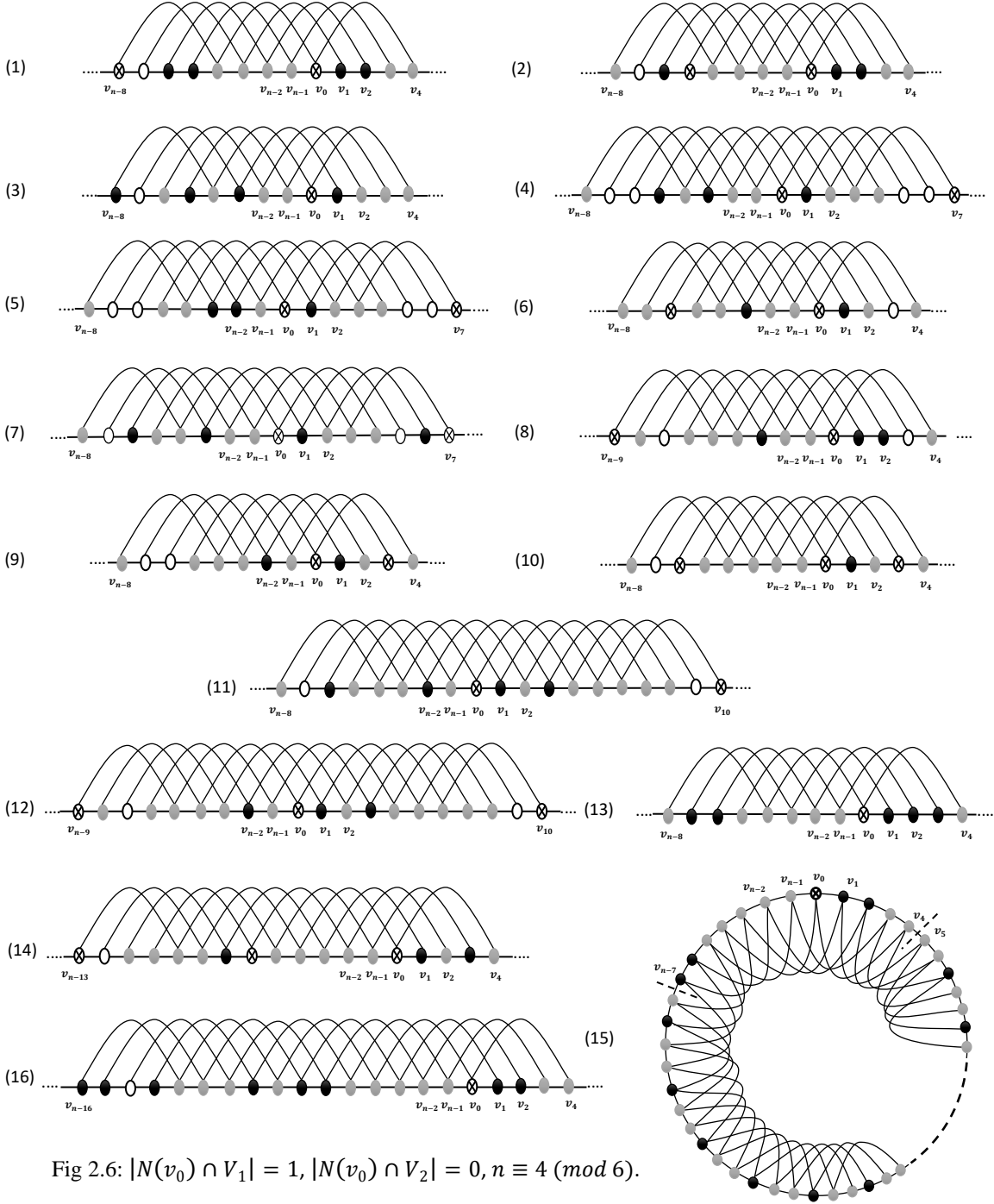

Fig 2.6: $|N(v_0) \cap V_1| = 1$, $|N(v_0) \cap V_2| = 0$, $n \equiv 4 \pmod{6}$.

From now on for every vertex $v$ from $V_2$ then $|N(v) \cap V_1| = 0$, and by continuing in the same way and without going through details we have:

**Case 4:** $v_0 \in V_2$ and $|N(v_0) \cap V_1| = 0$.

**Case 4.1:** $n \equiv 1, 2, 3, 5 \pmod{6}$, then $|V_2| \leq 2$, $|N(v_0) \cap V_2| = 0$ and $|V_{22}| = |V_{13}| = |V_{03}| = |V_{12}| = |V_{02}| = |V_{04}| = |E_2| = 0$, that leads to $v_1 \in V_{01} \cup V_{11} \cup V_{21} \cup V_{31}$ and we will discuss each one of them as follows:

- $v_1 \in V_{31}$ then $|V_{31}| = |V_2| = 1\sim$, but $v_6 \in V_{30} \cup V_{31} \cup V_{40}$ (see Fig 2.7 (1)), a contradiction.

- $v_1 \in V_{21}$, if $v_{n-3} \in V_0$ then $v_4 \in V_{11}$, $|V_{21}| = |V_2| = |V_{11}| = 1\sim$, and this leads to $v_{n-4}, v_{n-1} \in V_{01}$ and $v_{n-3}$ forces $v_{n-7} \in V_2$ (see Fig 2.7 (2)), a contradiction.

- $v_1 \in V_{11}$, if $v_2 \in V_1$ then $v_{n-2} \in V_0$ and $v_{n-3}$ forces $v_{n-7} \in V_2$ then $v_{n-6}, v_{n-8} \in V_0$, but $v_{n-2}$ leads to a contradiction with the 2RDF definition (see Fig 2.7 (3))

- $v_1 \in V_{01}$ and by symmetry $v_{n-1} \in V_{01}$, then $v_2, v_3$ force $v_6, v_7 \in V_2$ (see Fig 2.7 (4)), a contradiction.

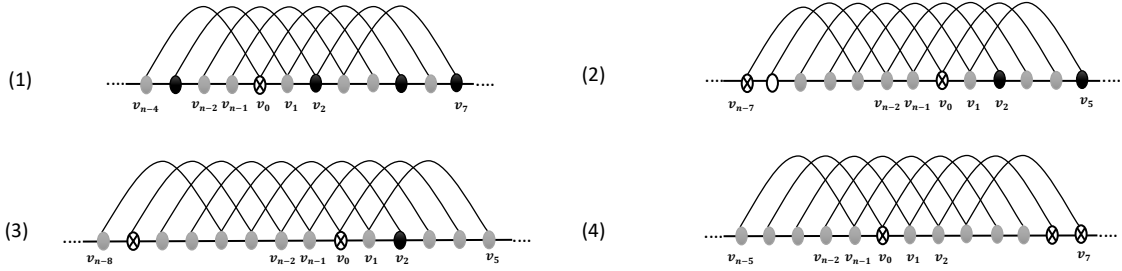

Fig 2.7: $|N(v_0) \cap V_1| = 0, n \equiv 1, 2, 3, 5 \pmod{6}$.

**Case 4.2:** $n \equiv 4 \pmod{6}$ then $|N(v_0) \cap V_2| \leq 1$.

**Case 4.2.1:** $|N(v_0) \cap V_2| = 1$, let $v_1 \in V_2$ then $v_0 v_1 \in E_2$ and $v_2, v_5, v_{n-3} \in V_0$ so we have the following cases:

- $v_{n-2} \in V_1$, then $v_{n-3}, v_{n-1}, v_2 \in V_{11}$ and that makes $|V_2| = 2$, $|E_2| = 1$, $|V_{11}| = 3\sim$, but $v_3$ forces $v_7 \in V_2$ (see Fig 2.8 (1)), a contradiction.

- $v_{n-2} \in V_0$, and this vertex makes $v_{n-6} \in V_2$, $|V_2| = 3$, $|E_2| = 1$, $|V_{30}| + |V_{11}| \leq 1\sim$, therefore $v_3 \in V_0$ and it forces $v_7 \in V_2$ (see Fig 2.8 (2)), a contradiction.

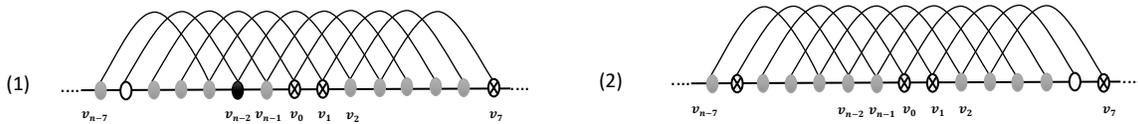

Fig 2.8: $|N(v_0) \cap V_1| = 0, |N(v_0) \cap V_2| = 1, n \equiv 4 \pmod{6}$.

**Case 4.2.2:** $|N(v_0) \cap V_2| = 0$, from previous discussion we conclude that for every vertex $v$ from $V_2$ then $N(v) \subset V_0$. (i)

That makes $v_1 \in V_{01} \cup V_{11} \cup V_{21} \cup V_{31} \cup V_{02} \cup V_{12} \cup V_{22}$ and will discuss them as follows:

- $v_1 \in V_{22}$, let $v_2 \in V_2$ then $v_7, v_8 \in V_0$, and $v_4 \in V_{11}$, and now $v_7$ forces $v_{11} \in V_2$ and that makes $|V_2| = 3$, $|V_{22}| = |V_{11}| = 1 \sim$, but $v_{n-4} \in V_{11} \cup V_{21} \cup V_{31} \cup V_{41}$ (see Fig 2.9 (1)), a contradiction.

- $v_1 \in V_{12}$, let $v_2 \in V_2$, $v_5 \in V_1$. Now, $v_{n-5} \in V_0$ which forces $v_{n-9} \in V_2$, and that makes $|V_2| = 4$, $|V_{12}| = 1 \sim$, but $v_{n-8} \in V_{02} \cup V_{12}$ (see Fig 2.9 (2)), a contradiction.

- $v_1 \in V_{02}$, let $v_2 \in V_2$ then $v_{n-3}, v_5$ force $v_{n-7}, v_9 \in V_2$ and that makes $|V_2| = 4$, $|V_{02}| = 1$, $|V_{30}| + |V_{11}| \leq 1 \sim$, but $v_7 \in V_0$ force $v_{11} \in V_2$ (see Fig 2.9 (3)), a contradiction.

- $v_1 \in V_{31}$, $v_{n-1} \in V_{01}$, $v_{n-8}, v_8 \in V_0$ then $v_{n-4}, v_4 \in V_{11}$ and that makes $v_7 \in V_1$ and $v_6 \in V_0$, moreover $v_5$ makes $v_{n-6}, v_{n-9} \in V_1$ and $v_{n-2} \in V_{30}$, and by continuing in this way for $8 < j < n$ we have
$$\begin{cases} v_j \in V_0: j = 6i, 6i+2, 6i+3, 6i+5. \\ v_j \in V_1: j = 6i+1, 6i+4. \end{cases}$$

And here $v_5$ leads to a contradiction with 2RDF definition (see Fig 2.9 (4)), a contradiction.

- $v_1 \in V_{21}$, let $v_5 \in V_0$ then $v_{n-1} \in V_{01}$ that makes the vertices $v_3$, $v_5$ force $v_7, v_6, v_9 \in V_1$, thus $v_6 v_7 \in E_1$ and $v_{n-5}$ makes $v_{n-6}, v_{n-9} \in V_1$, so $v_{n-2} \in V_{30}$, and by continuing in this way for $9 < j < n$ we have

$$\begin{cases} v_j \in V_0: j = 6i, 6i+2, 6i+3, 6i+5, \\ v_j \in V_1: j = 6i+1, 6i+4. \end{cases}$$

And that makes $v_9 v_{10}, v_9 v_{13}, v_6 v_{10} \in E_1$ then $2|V_2| + 2|V_{21}| + 2|E_1| + |V_{11}| + |V_{30}| > 11$ (see Fig 2.9 (5)), a contradiction.

- $v_1 \in V_{11}$, let $v_2 \in V_1$ then $v_{n-1} \in V_{01}$ and $v_{n-2}$ forces $v_{n-6} \in V_1$, implies that $v_{n-3}$ makes $v_{n-7} \in V_2$, a contradiction with (i) (see Fig 2.9 (6)).

- $v_1 \in V_{01}$ then $v_{n-1} \in V_{01}$ and $v_2, v_3$ force $v_6, v_7 \in V_2$ then $v_6 v_7 \in E_2$, a contradiction with (i) (see Fig 2.9 (7)).

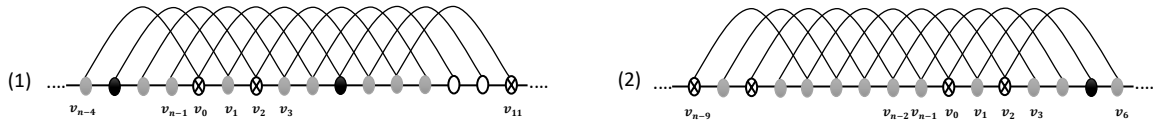

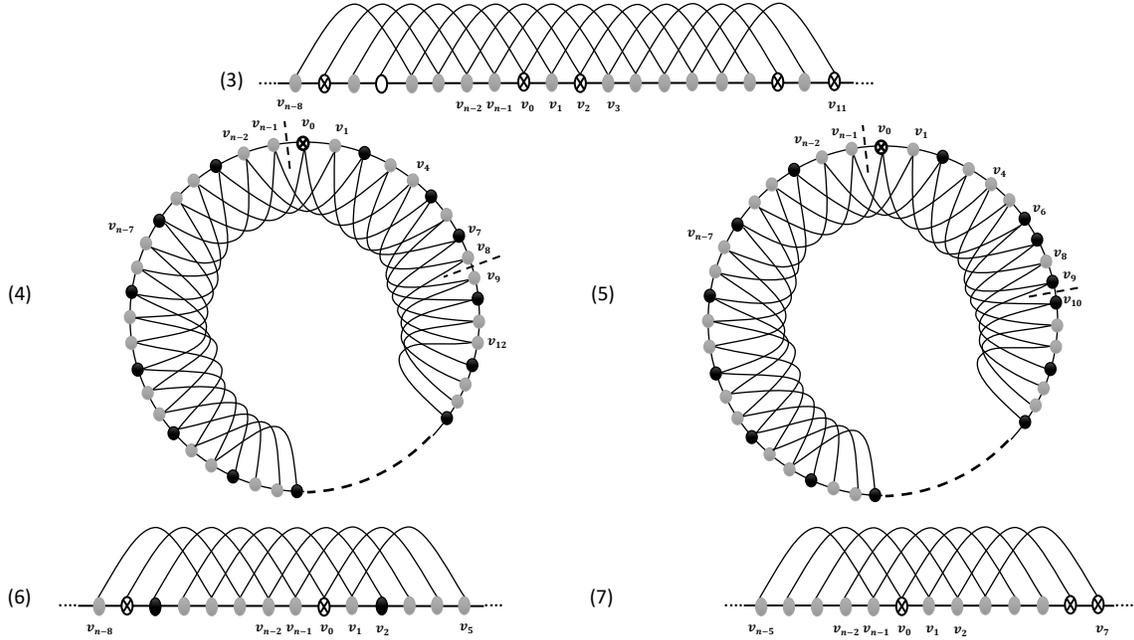

Fig 2.9: $|N(v_0) \cap V_1| = 0, |N(v_0) \cap V_2| = 0, n \equiv 4 \pmod 6$.

**Case 5:** $|V_2| = 0$.
Then $|V_{11}| = |V_{12}| = |V_{13}| = |V_{21}| = |V_{22}| = |V_{31}| = |V_{02}| = |V_{03}| = |V_{04}| = |E_{12}| = |E_2| = 0$, and $\beta = 2|E_1| + 2|V_{40}| + |V_{30}|$.

**Case 5.1:** $n \equiv 1, 2, 3, 5 \pmod 6$, then $|E_1| \le 2$ and $|V_{40}| \le 2$.

**Case 5.1.1:** $v_0 \in V_{40}$, so if $v_1 v_2 \in E_1$ then $|V_{40}| = 1$, $|E_1| = 1$ and $|V_{30}| \le 1$ which leads to $v_{n-8}, v_{n-5}, v_{n-3}, v_{n-2}, v_3, v_5, v_6, v_7, v_8 \in V_0$ and $|V_{30}| = 1$, but $v_7$ forces $v_{11} \in V_2$ (see Fig 2.10 (1)), a contradiction. Then $v_2 \in V_0$ and by symmetry we have $v_{n-2} \in V_0$. Now, $v_{n-5}, v_{n-3}, v_3, v_5 \in V_0$, if $v_9 \in V_1$ then $v_{n-5} \in V_{30}$ and $v_5 \in V_{40}$, moreover $|V_{40}| = 2$, $|V_{30}| = 1$ and that makes $v_{n-7}, v_{n-8}, v_{n-9}, v_{n-10} \in V_0$, but here $v_{n-9}$ forces $v_{n-13} \in V_2$ (see Fig 2.10 (2)), a contradiction, from that we conclude that $v_9 \in V_0$ and by symmetry we have $v_{n-9} \in V_0$ and that makes $v_5, v_{n-5} \in V_{30}$, $|V_{40}| = 1$, $|E_1| = 0$ and $|V_{30}| \le 3$, but $v_{n-9}$ forces $v_{13} \in V_2$ (see Fig 2.10 (3)), a contradiction.

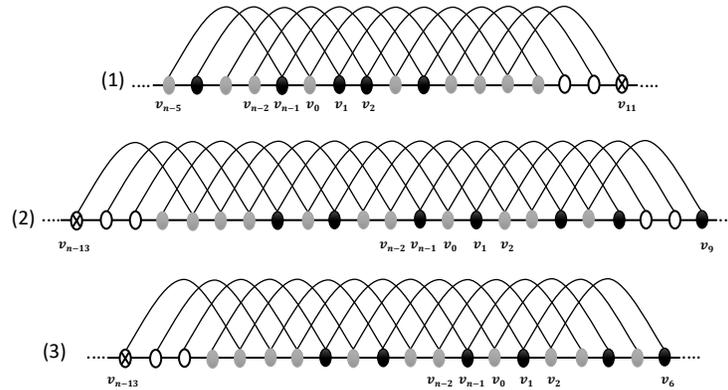

Fig 2.10: $|V_2| = 0, n \equiv 1, 2, 3, 5 \pmod 6, v_0 \in V_{40}$.

From now on $|V_{40}| = 0$ and $\beta = 2|E_1| + |V_{30}|$ in Case 5.1.

**Case 5.1.2:** $v_0 v_1 \in E_1$, Therefore we will discuss each of the following vertices:

- $v_2 \in V_1$, then $|E_1| = 2$, $|V_{30}| \leq 1$ and that leads to $v_3, v_4, v_5, v_6 \in V_0$, so $v_4, v_5$ force $v_8, v_9 \in V_1$ which makes $v_8 v_9 \in E_1$ (see Fig 2.11 (1)), a contradiction. Thus, we have $v_2 \in V_0$ and we conclude that for every vertex $v_i$ from $V(G)$ which implies that $v_i v_{i+1} \in E_1$ makes $v_{i-1}, v_{i+2} \in V_0$, and by that we have $v_{n-1} \in V_0$.

- $v_5 \in V_1$ then $|E_1| = 2$, $|V_{30}| \leq 1$, and that makes $v_{n-4}, v_{n-3}, v_4, v_6, v_9 \in V_0$, so $v_2$ forces either $v_3 \in V_1$ then $v_7, v_8 \in V_0$, $v_4 \in V_{30}$ and $|V_{30}| = 1$, but $v_8$ makes $v_{12} \in V_2$ (see Fig 2.11 (2)), a contradiction, or $v_{n-2} \in V_1$ then $v_3$ forces $v_7 \in V_2$ a contradiction, then $v_5 \in V_0$ and by symmetry we have $v_{n-4} \in V_0$.

- $v_4 \in V_1$ then $|E_1| = 2$, $|V_{30}| \leq 1$, that leads to $v_7 \in V_1$, $v_{n-5}, v_{n-3}, v_{n-6}, v_3, v_6, v_8 \in V_0$ and here $v_{n-4}$ force $v_{n-8} \in V_1$ which makes $v_{n-12}, v_{n-9}, v_{n-7} \in V_0$, but $v_{n-5}$ leads to a contradiction with the 2RDF definition (see Fig 2.11 (3)), so we have $v_4 \in V_0$ and by symmetry $v_{n-3} \in V_0$, therefore $v_{n-2}, v_3 \in V_1$, $v_2, v_{n-1} \in V_{30}$, $|E_1| = 1$ and $|V_{30}| \leq 3$ thus we have $v_{n-6}, v_7 \in V_0$ and $v_6$ forces $v_{10} \in V_2$ (see Fig 2.11 (4)), a contradiction. From now on $|E_1| = 0$ in Case 5.1.

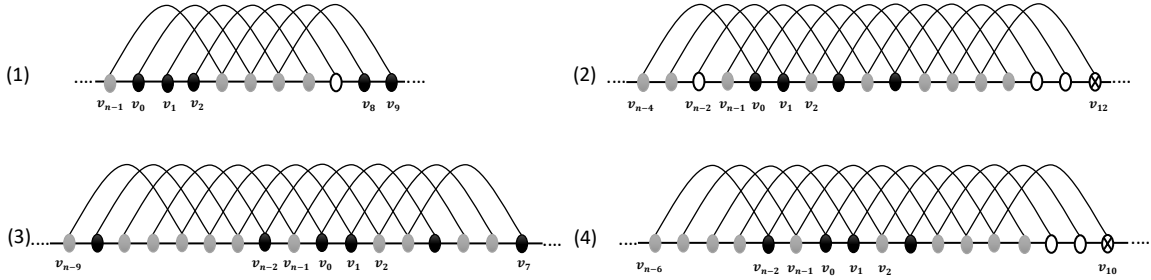

Fig 2.11: $|V_2| = 0$, $n \equiv 1, 2, 3, 5 \pmod 6$, $v_0 v_1 \in E_1$.

**Case 5.1.3:** $v_0 \in V_{30}$, let $v_{n-1} \in V_0$ then $v_{n-8}, v_{n-5} v_{n-3}, v_2, v_3, v_5, v_8 \in V_0$, but $v_{n-1}$ makes $v_{n-2} \in V_2$ (see Fig 2.12), a contradiction.

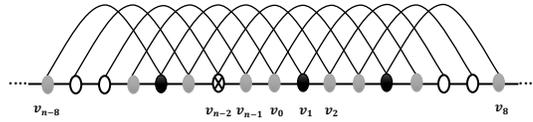

Fig 2.12: $|V_2| = 0$, $n \equiv 1, 2, 3, 5 \pmod 6$, $v_0 \in$

**Case 5.1.4:** $|V_{30}| = 0$, $v_0 \in V_0$ If $v_1, v_{n-1} \in V_1$ then $v_{n-5}, v_{n-3}, v_{n-2}, v_2, v_3, v_5 \in V_0$, therefore $v_4$ makes $v_8 \in V_2$ a contradiction. And if $v_1, v_4 \in V_1$ then $v_{n-4}, v_{n-1}, v_2, v_3, v_5, v_8 \in V_0$ and here $v_3$ forces $v_7 \in V_1$ and that leads to $v_4, v_{10} \in V_0$, the vertices $v_2, v_6$ makes $v_{n-2}, v_{10} \in V_1$ and by continuing in this way for $0 < j < n - 4$ we have
$$\begin{cases} v_j \in V_0: j = 6i, 6i + 2, 6i + 3, 6i + 5, \\ v_j \in V_1: j = 6i + 1, 6i + 4 \text{ and } d(v_{6i+1}) = 1, d(v_{6i+4}) = 2. \end{cases}$$
And when:

- $n \equiv 1, 2, 5 \pmod{6}$, then $v_{n-4}$ leads to a contradiction with the 2RDF definition (see Fig 2.13).
- $n \equiv 3 \pmod{6}$, then $v_{n-6}$ leads to a contradiction with the 2RDF definition (see Fig 2.13).

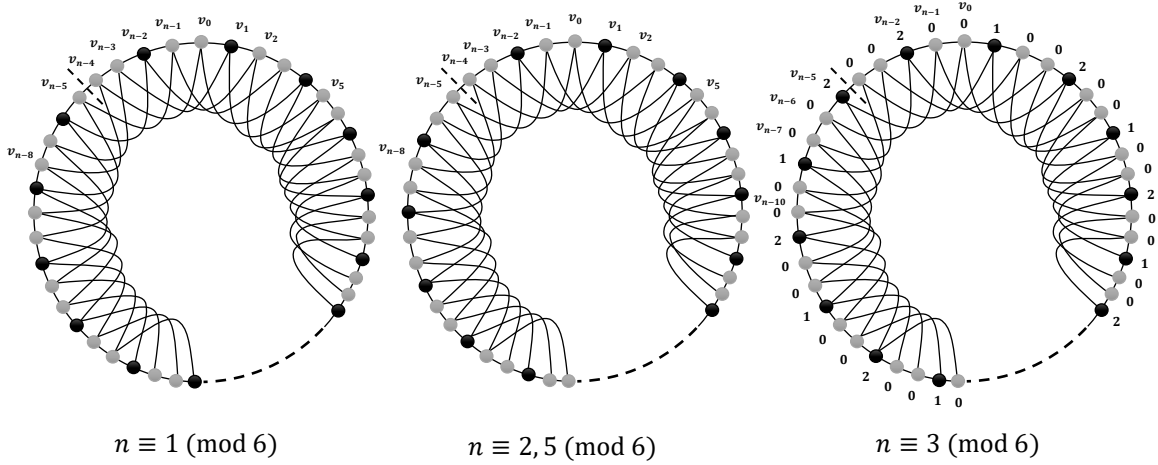

$n \equiv 1 \pmod{6}$      $n \equiv 2, 5 \pmod{6}$      $n \equiv 3 \pmod{6}$

Fig 2.13: $|V_2| = |V_{40}| = |V_{30}| = 0, n \equiv 1, 2, 3, 5 \pmod{6}$.

**Case 5.2:** $n \equiv 4 \pmod{6}$, let $v_0 \in V_1$.

**Case 5.2.1:** $|N(v_0) \cap V_1| = 4$ then $v_2, v_3, v_9 \in V_0$, $|E_1| = 4$, $|V_{40}| = 0$ and $|V_{30}| \le 3$, that leads to $v_{n-10}, v_{n-8}, v_{n-7}, v_7, v_8, v_{10} \in V_0$, but $v_9$ forces $v_{13} \in V_2$ (see Fig 2.14 (1)), a contradiction.

**Case 5.2.2:** $|N(v_0) \cap V_1| = 3$, let $v_1 \in V_0$, then $v_{n-5}, v_{n-3}, v_{n-2}, v_3, v_5, v_6, v_7, v_8 \in V_0$ and $v_2, v_9, v_{10} \in V_1$, but $v_7$ forces $v_{11} \in V_2$ (see Fig 2.14 (2)), a contradiction.

**Case 5.2.3:** $|N(v_0) \cap V_1| = 2$, let $v_1, v_{n-1} \in V_1$ then $v_{n-3}, v_{n-2}, v_2, v_3 \in V_0$ and that makes $v_{n-7}, v_{n-6}, v_6, v_7 \in V_1$ and $v_{n-5}, v_5 \in V_0$, but $v_{n-4}, v_4$ force $v_{n-8}, v_8 \in V_1$ then $2|E_1| > 11$ (see Fig 2.14 (3)), a contradiction.

**Case 5.2.4:** $|N(v_0) \cap V_1| = 1$, let $v_1 \in V_1$ then $v_{n-2}, v_{n-1}, v_6, v_{n-5}, v_9, v_{12} \in V_1$ $v_{n-6}, v_7, v_8 \in V_0$, $v_{n-1}, v_2 \in V_{40}$, $v_5 \in V_{30}$ and by continuing in the same way for $3 < j < n-6$ we have
$$\begin{cases} v_j \in V_0 : j = 6i+1, 6i+2, 6i+4, 6i+5, \\ v_j \in V_1 : j = 6i, 6i+3. \end{cases}$$
And $v_{n-8}$ leads to a contradiction with the 2RDF definition (see Fig 2.14 (4)).

**Case 5.2.5:** $|N(v_0) \cap V_1| = 0$ then $v_{n-2}, v_2 \in V_0$ and $v_2$ forces $v_3, v_6 \in V_1$, continuing in the same way for $0 < j < n-4$ we have
$$\begin{cases} v_j \in V_0 : j = 6i+1, 6i+2, 6i+4, 6i+5, \\ v_j \in V_1 : j = 6i, 6i+3. \end{cases}$$

And $v_{n-5}$ leads to a contradiction with the 2RDF definition (see Fig 2.14 (5)).

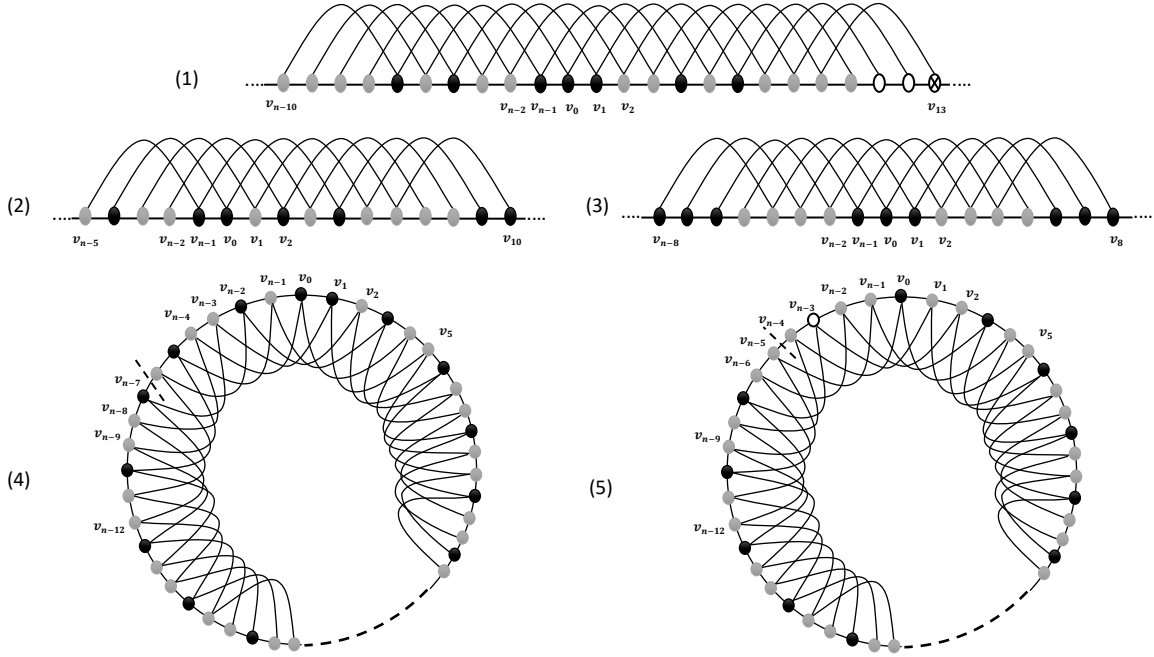

Fig 2.14: $|V_2| = 0$, $n \equiv 4 \pmod 6$.

And by that discussion we see that every possible Case leads to a contradiction with our assumption. □

**Theorem 2.4.**

For $n \geq 9$ then

$$\gamma_{r2}\left(C(n; \{1,4\})\right) = \begin{cases} \left\lceil \dfrac{n}{3} \right\rceil: & n \equiv 0 \pmod 6, \\ \left\lceil \dfrac{n}{3} \right\rceil + 1: & n \equiv 1,2,3,5 \pmod 6, \\ \left\lceil \dfrac{n}{3} \right\rceil + 2: & n \equiv 4 \pmod 6. \end{cases}$$

**proof.** We found by Lemma 2.1. that

$$\gamma_{r2}\left(C(n; \{1,4\})\right) \leq \begin{cases} \left\lceil \dfrac{n}{3} \right\rceil: & n \equiv 0 \pmod 6, \\ \left\lceil \dfrac{n}{3} \right\rceil + 1: & n \equiv 1,2,3,5 \pmod 6, \\ \left\lceil \dfrac{n}{3} \right\rceil + 2: & n \equiv 4 \pmod 6. \end{cases}$$

And we conclude by Theorem 3. that $\gamma_{r2}(C(n;\{1,4\})) \geq \left\lceil \frac{n}{3} \right\rceil$ when $n \equiv 0 \pmod{6}$, also by Lemma 2.4. we found that $\gamma_{r2}(C(n;\{1,4\})) \geq \left\lceil \frac{2n+\beta}{6} \right\rceil \geq \left\lceil \frac{2n+6}{6} \right\rceil = \left\lceil \frac{n}{3} \right\rceil + 1$ when $n \equiv 1,2,3,5 \pmod{6}$ and $\gamma_{r2}(C(n;\{1,4\})) \geq \left\lceil \frac{2n+\beta}{6} \right\rceil \geq \left\lceil \frac{2n+12}{6} \right\rceil = \left\lceil \frac{n}{3} \right\rceil + 2$ when $n \equiv 4 \pmod{6}$. □

**Acknowledgements**

**References**

[1] B. Brešar, M. A. Henning and D. F. Rall, Paired-domination of Cartesian products of graphs and rainbow domination. Electron. Notes Discrete Math. 22 (2005), 233-237.
[2] B. Brešar and T. K. Sumenjak, On the 2-rainbow domination in graphs, Discrete Appl. Math. 155 (2007), 2394-2400.
[3] T. Chunling, L. Xiaohui, Y. Yuansheng and L. Meiqin. 2-rainbow domination of generalized Petersen graphs $P(n,2)$. Discrete Appl. Math. 157 (2009), 1932-1937.
[4] R. Erveš, J. Žerovnik, On 2-Rainbow Domination Number of Generalized Petersen Graphs $P(5k,k)$. Symmetry. 13, 809 (2021).
[5] B. Brešar, M. A. Henning and D. F. Rall. Rainbow domination in graphs. Taiwanese J. Math. 12 (2008), 213-225.
[6] G. Xu. 2-rainbow domination in generalized Petersen graphs $P(n,3)$. Discrete Appl. Math. 157 (2009) 2570-2573.
[7] T. K. Sumenjak, D. F. Rall and A. Tepeh, Rainbow domination in the lexicographic product of graphs, Discrete Appl. Math. 161 (2013), 2133–2141.
[8] Z. Stepień and M. Zwierzchowski, 2-rainbow domination number of Cartesian products: Cn □ C3 and Cn □ C5, J. Comb. Optim. 28 (2014), 748–755.
[9] Z. Stepień, A. Szymaszkiewicz and L. Szymaszkiewicz, 2-Rainbow domination number of Cn □ C5, Discrete Appl. Math. 170 (2014), 113–116.
[10] Z. Shao, H. Jiang, P. Wub, S. Wangd, J. Žerovnik, X. Zhang and J. Liu, On 2-rainbow domination of generalized Petersen graphs, Discrete Appl. Math. accepted, (2018).
[11] D. A. Mojdeh and Z. Mansouri, Rainbow domination of graphs, International Conference on Combinatorics, Cryptography and Computation. I4C, (2018).
[12] M. Ali, M. T. Rahim, M. Zeb and G. Ali, On 2-rainbow domination of some families of graphs. International Journal of Mathematics and Soft Computing 1, (2011), 47–53.
[13] Z. Shao, M. Liang, C. Yin, X. Xu, P. Pavlič and J. Žerovnik, On rainbow domination numbers of graphs. Inform. Sci. (2014), 254, 225–234.
[14] X. Fu, X. Wu, G. Dong, H. Li and W Guo, 2-Rainbow Domination of the Circulant Graph $C(n;\{1,3\})$, International Conference on Test, Measurement and Computational Method, TMCM (2015).

**[15]** N. Obradović, J. Peters and G. Ružić, Efficient domination in circulant graph with two chord lengths, Inform. Process. Lett. 102 (2007) 253–258.

**[16]** K. Kumar and G. MacGillivray, Efficient domination in circulant graphs, Discrete Math 313 (2013) 767–771.

**[17]** H. A. Ahangar, J. Amjadi, S. M. Sheikholeslami and D. Kuziak, Maximal 2-rainbow domination number of a graph, AKCE Int. J. Graphs and combin. 13 (2016), 157–164.

**[18]** J. Amjadi, S. M. Sheikholeslami and L. Volkmann, Rainbow restrained domination numbers in graphs, Ars Combin., 124 (2016), 3–19.

**[19]** Y. Rao, P. Wu, Z. Shao, R. Shaheen, S.M. Sheikholeslami and L.Chen, The Rainbow Restrained Domination in Torus Network, International Conference on Cyber-Enabled Distributed Computing and Knowledge Discovery, (2018).

**[20]** H. A. Ahangar, J. Amjadi, M. Chellali, S. Nazari-Moghaddam, S.M. Sheikholeslami, Total 2-Rainbow Domination Numbers of Trees. Discuss. Math. Graph Theory 41 (2021), 345–360.

**[21]** H. Jiang, Y. Rao, Total 2-Rainbow Domination in Graphs, Mathematics, 10 (2022), 2059.

**[22]** Z. Mansouri and D.A. Mojdeh, Outer independent rainbow dominating functions in graphs, Opuscula Math. 40(5)(2020), 599-615.

**[23]** N. J. Rad, E. Gholami, A. Tehranian, H. Rasouli, A new upper bound on the independent 2-rainbow domination number in trees, Commun. Combin. Optim. (2022).